\documentclass[3p, 10pt , final]{elsarticle}
\usepackage{graphicx}
\usepackage{amssymb}
\usepackage{hyperref}
\usepackage{amsmath}

\newcommand{\Rb}{\mathbb{R}}

\newtheorem{proposition}{Proposition}[section] 
\newtheorem{definition}{Definition}[section]

\newtheorem{lemma}{Lemma}[section]
\newtheorem{thm}{Theorem}[section]
\newtheorem{exemplo}{Example}[section]
\newtheorem{remark}{Remark}[section]
\begin{document}
\begin{frontmatter}
 \title{Root Finding by High Order Iterative Methods Based on Quadratures}
\author[b]{M\'{a}rio M. Gra\c{c}a}
\address[b] {LAETA,IDMEC, Departamento de Matem‡tica,  Instituto Superior T\'{e}cnico,  Universidade de Lisboa,
  1049-001 Lisboa,  Portugal}
  \ead{mgraca@math.ist.utl.pt}
  
\author[3]{Pedro M. Lima \corref{cor1}}
\address[3]{CEMAT,  Departamento de Matem‡tica, Instituto Superior T\'{e}cnico,  Universidade de Lisboa,
  1049-001 Lisboa,  Portugal}
\cortext[cor1]{Corresponding  author}
\ead{plima@math.ist.utl.pt}
 
\begin{abstract}
We discuss a recursive family of iterative methods for the numerical approximation of roots of nonlinear functions in one variable.
These methods are based on Newton-Cotes closed quadrature rules.
We prove that when a quadrature rule with
$n+1$ nodes is used  the resulting iterative method has convergence order at least $n+2$, starting with the case $n=0$ (which corresponds to the Newton's method). 
\end{abstract}
\end{frontmatter}
\section{Introduction}\label{introd}
The use of quadrature rules for the construction of iterative methods,
applied to the solution of nonlinear equations or systems, has been 
considered by many authors (see, for example, \cite{fernando}, \cite{frontini}, \cite{hafiz}, \cite{mir}).
However, so far these methods  have in general been treated separately or dealing with a specific quadrature rule or a small set of them. In this work we shall treat this matter in a systematic and unifying way. 

Our main purpose is to obtain a family of recursive iterative methods based on quadratures, with higher
convergence order than the Newton's method. This one is universally considered the method of choice for approximating a root  $z\in \Rb^d$ for a given  equation   $f(x)=0$, where $f: D\subset \Rb^d \mapsto \Rb^d$.
However, its limitations are also known. In many cases of practical interest, the Newton's method  fails to converge unless the initial approximation  lies in a small neighborhood of the root we want to approximate.  It's mainly in  such cases that higher order iterative methods can be useful, such as those described in the present paper.

For a fixed positive integer $n$,   we define recursively a certain function $t_n: \Rb \mapsto \Rb$, based on a Newton-Cotes closed quadrature rule, with $n+1$ nodes (see Definition  \ref{def1}).   Numerical integration is discussed, for example, in \cite{krylov} Ch. 6, \cite{brass} Ch. 5, and for  Newton\--Cotes quadrature rules,  see for instance \cite{gautschi},  Ch. 3,  and \cite{graca1}. 

In the present work, we take as the basic iteration function $t_0$
the Newton's iterative process.
Assuming that the Newton's method has convergence order $p \ge 2$, we prove in Theorem \ref{principal} that the convergence order of  our iterative function $t_n$ is not less than $2+n$. This enables us to 
construct iterative methods of arbitrary convergence order 
for the numerical solution of nonlinear equations. 

In Sec. \ref{rel}   we establish a relation between  the  approximation of a real root $z$ of an equation
$f(x)=0$ and a quadrature rule, using the main theorem of integral calculus. Though   we restrict ourselves to closed
Newton\--Cotes rules,   quadratures of different types can be also used.

The iterative methods described here can be easily extended to the case of multivariate 
functions. However, the analysis of convergence in this case is out of the scope of the present paper.

In Sec. \ref{sec new} we begin by showing that the iterative function $t_0$
 coincides with the classical Newton's iterative function. As known,
if $z$ is a simple root, this  process has at least second order convergence, provided that the initial approximation is sufficiently close to a simple root $z$. 
For a given integer $n\geq 1$, we show how to apply a certain closed  Newton\--Cotes quadrature rule (with $n+1$ nodes), in order   that the corresponding  recursive iterative method, $t_n$, possesses in general a higher convergence order than the previous iterative mapping $t_{n-1}$.
 Namely, we prove that each iterative mapping $t_n$ has convergence order
 not less than $n+2$, which is the main result of this work.\footnote{In certain particular cases the order of 
$t_n$ may be higher than $n+2$; in such cases it may happen that $t_{n}$ and  $t_{n+1}$ have the same convergence order $n+3$.}

In Sec.~4  we present some numerical examples illustrating the application of the described methods.  We compare their accuracy  and verify experimentally their convergence order. Special attention is paid to the cases where the classical Newton's method fails.

Finally, in Sec.~5, we present the main conclusions and discuss perspectives for a future work.
   
\section{Iterative Methods for Root Finding and Quadrature Rules}\label{rel}

Given a  function $f$ in one real variable, let $z$ be a simple  root  of $f$ (that is $f'(z) \ne 0$). Suppose that $f$ is sufficiently regular in a certain neighborhood of $z$.
By the fundamental theorem of integral calculus we know that

\begin{equation}\label{q1}
\int_{x}^ z f^ {(1)} (t) \, dt=f(z)-f(x)= - f(x).
\end{equation}
Choosing a non\--negative integer $n$,
we approximate the integral on the left-hand side of \eqref{q1}  by a certain interpolatory quadrature rule with $n+1$ nodes, which we denote by $Q_n\left(f^ {(1)}\right)$. We can write the rule as
 \begin{equation}\label{q2}
Q_n\left(f^ {(1)}\right) = \displaystyle{ \frac{z-x}{c_n}  } \, B_n(x).
 \end{equation}
 The function  $B_n(x)$  in  \eqref{q2}
 is defined by the weights $A_i$ and by the \lq\lq{nodes}\rq\rq $\,\xi_i(x)\in [x,z]$, such that
  \begin{equation}\label{q3}
 B_n(x)= A_0\, f^ {(1)} (\xi_0(x))+ A_1\, f^ {(1)} (\xi_1(x))+\cdots+  A_n\, f^ {(1)} (\xi_n(x)).
 \end{equation}
In an interpolatory quadrature rule, the constant
 $c_n$ in \eqref{q2} satisfies the equality
 \begin{equation}\label{q4}
 c_n=\sum_{i=0}^ n A_i,
 \end{equation}
 since, by construction, the rule is exact when applied to
  $f(x)\equiv 1$. \footnote{ We assume that the length of the interval where the quadrature rule is applied
is $c_n$, where $c_n$ is the least integer for which all the weights  $A_i$ are integer numbers. When we consider integration on an interval of a different length, all the weights should be multiplied by a certain number, explaining why the factor
 $(z-x)/c_n$ appears in formula (\ref{q2}).}
 
 \medskip
 \noindent
We also assume that for $t\in [x,z]$, the function $B_n(t)$ is finite and has a finite inverse, that is,
   \begin{equation}\label{q5}
B_n(t)\neq 0 \quad \forall t\in [x,z].
 \end{equation}
Finally, the quadrature nodes
 $\xi_i(x)$ satisfy  
 \begin{equation}\label{q10}
 \xi_i(z)=z,\quad \mbox{para}\quad i=0,1, \ldots, n.
 \end{equation}
 In Sec. \ref{sec new} we will define the functions
 $\xi_i(x)$, which are the quadrature nodes  in \eqref{q3}, 
 using the  closed   Newton\--Cotes quadrature rules (\cite{gautschi}, Ch. 3).
 
 \medskip
 The iterative processes to be constructed will possess some of the properties of the adopted quadrature rules, and this will be reflected in the following proofs. In a future work we intend to use open quadrature rules with the same purpose. 
  
 \medskip
 \noindent
  Substituting \eqref{q2} into \eqref{q1}, we obtain
 \begin{equation}\label{q6}
 z-x\simeq -c_n\, B_n^ {-1} (x)\, f(x).
 \end{equation}
The approximate equality \eqref{q6} leads us to the following definition of the mapping $t_n$.

 \begin{definition}(Iterative mapping based on a quadrature rule)\label{def1}
 
For a given integer $n\geq 0$ and a certain  function  $B_n(x)$, associated to the quadrature rule  \eqref{q2}, satisfying the conditions \eqref{q3}\--\eqref{q10}, 
the iterative mapping $t_n$
is  defined by
 \begin{equation}\label{q7}
 t_n(x)=x - c_n \, B_n^ {-1} (x) f(x).
 \end{equation}
 \end{definition}
Defining  the auxiliary function
\begin{equation}\label{q8}
H_n(x)= t_n(x)-x,
\end{equation}
we remark that $H_n$ satisfies
\begin{equation}\label{q9}
B_n(x)\, H_n(x)=- c_n \, f(x) \Longleftrightarrow H_n(x)= -c_n\, B_n^ {-1}(x)\, f(x).
\end{equation}

 \noindent
 Since, for $n\geq 1$, we will use only closed Newton\--Cotes
quadrature rules, the  function 
 $t_n$ in \eqref{q7} will be called the {\em  Newton-Cotes closed iterative mapping with $n+1$ nodes}.

\medskip
\noindent 
We begin by proving the  superlinear convergence of the  mapping defined by \eqref{q7}, in the case 
$f$ is a  one-variable function, sufficiently regular in the neighborhood of a simple root $z$. 
 
 \begin{proposition}\label{conv} (Superlinear convergence of iterating mappings)

 \noindent
 A simple root of the equation $f(x)=0$ is a fixed point
 of the iterative mapping \eqref{q7}. Moreover, starting
from an approximation  $x_0$ sufficiently close to $z$, the sequence defined by $x_{k+1}= t_n(x_k)$
 converges superlinearly to $z$, for any $n\geq 0$.
 \end{proposition}

{ \bf Proof.}
From \eqref{q3}, taking the equalities  \eqref{q10} into account,  we obtain
 $$
 \begin{array}{ll}
 B_n(z)&= A_0\, f^ {(1)} (z)+ A_1\, f^ {(1)} (z)+ \cdots+  A_n\, f^ {(1)} (z)\\
 &=  f^ {(1)} (z)\, \sum_{i=0}^ n A_i.
 \end{array}
 $$
Since, by construction, the sum of the weights $A_i$ is equal to
 $c_n$,  it follows that
 \begin{equation}\label{q11}
 B_n(z)= c_n\, f^ {(1)} (z),
 \end{equation}
 and therefore $B_n(z)\neq 0$, since $ z$ is a simple root of $f$.
Moreover, from  \eqref{q7},  we have
 $$
 t_n(z)= z- c_n\, c_n^ {-1} \left( f^ {(1)} (z) \right)^ {-1} f(z)=z,
 $$
which means that a simple root of $f$ is a fixed point of $t_n$.
From \eqref{q8}, we then conclude that $H_n$ vanishes at the fixed point $z$:
 \begin{equation}\label{q12}
 H_n(z)=0.
 \end{equation}
Differentiating both sides of \eqref{q9}, we obtain
 \begin{equation}\label{q13}
 B_n^ {(1)} (x)\, H_n(x)+ B_n(x)\, H_{n}^ {(1)} (x)= - c_n\, f^ {(1)} (x).
 \end{equation}
Hence, taking   \eqref{q12} into consideration, 
 from \eqref{q13} we conclude that
 $$
 B_n(z)\, H_n^ {(1)} (z)= -c_n\, f^ {(1)} (z).
 $$
 From the last equality, knowing that  $B_n(z)$ satisfies \eqref{q11}, we get
 $$
 c_n\, f^ {(1)} (z)\, H_n^ {(1)} (z)= -c_n\, f^ {(1)} (z),
 $$
 or, taking  \eqref{q8} into consideration,
 $$
  H_n^ {(1)} (z)=-1\quad  \Leftrightarrow \quad t_n^ {(1)} (z)-1=-1\quad \Leftrightarrow \quad t_n^ {(1)}(z)=0.
 $$
The last equality means that the iterative process generated
by  $t_n$ converges locally to $z$ and the convergence is superlinear.
 $\hfill \Box$
 
 \medskip

  Once an iterating mapping $t_0$ is chosen, having
superlinear convergence,  the Proposition  \ref{conv}  enables us to construct other mappings, based on quadrature rules, whose convergence order is not less than 2 (the same convergence order as the Newton's method, when applied to a simple root).
Moreover, by an adequate choice of the nodes of the quadrature rule
 $B_n$,  following Definition  \ref{def1}, we can build new methods whose convergence order is higher than 2.

By modifying the function $f$ (as described in the next subsection), we can also deal with the case of a multiple root. 
Therefore recursive iterative mappings $t_n$ can be obtained, having an arbitrarily high order, provided the
mapping  $t_0$ is chosen  so that it converges superlinearly to the considered root $z$.
  
\subsection{Multiple Roots}\label{secmult}

It is a common technique to modify a given function $f$ if the Newton's method does not provide satisfactory results, when applied to  
its roots (see, for example, \cite{benisrael1}, \cite{graca} and references therein). For example, if $t_0$  is the Newton's iterative mapping, for a function  $f$ with a multiple root $z$, 
one can define
$$
F(x)=t_0(x)-x=-\displaystyle{ \frac{f(x)}{f^ {(1)}(x)} }.
$$
Then  if $f''(z) \ne 0$ it is easy to show that $z$ is a simple root of $F$.
Therefore, Proposition  \ref{conv} holds in the case of multiple roots, provided that we start with the Newton's iterative mapping  $t_0$
applied to  $F$  (instead of the original function $f$) (see Example \ref{exemplo2}).
\section{Newton's, Trapezoidal and Simpson's Rules}\label{sec new}
\noindent
In this section we  introduce iterative functions
  $t_0$, $t_1$ and $t_2$, in $\Rb$, based on well-known quadrature rules. 
  The first of these functions results immediately from the application of the
  left rectangles rule (the only open Newton-Cotes rule considered in this paper);
the second one follows from $t_0$ and from the trapezoidal rule; finally the function $t_2$
results from $t_1$ and the Simpson's rule.  
Note that once $t_0$ has convergence order $p \ge 2$, the maps $t_1$ and  $t_2$ will have, by construction,   convergence orders at least 3 and 4, respectively.
 
 \medskip
 \noindent
 In Table \ref{tabnc} the weights  $A_i$  and the constants $c_n$ are displayed, needed for the construction of the  Newton\--Cotes iterative functions $t_n$, with $0 \leq n \leq 7$. We do not consider the case $n\geq 8$, since the weights $A_i$ may become negative for such values of $n$, which leads to numerically unstable formulae (see, for example,  \cite{bjorck1}, p. 534).
    \begin{table}
$$
\begin{array}{| c || c | c|  c|  c|  c| c|  c| c|| c| }
\hline
n&A_0& A_1& A_2& A_3& A_4& A_5& A_6& A_7 &c_n=\sum_{i=0}^ n A_i\\
\hline
0& 1 &  &  &  &  &  &  &  &   1\\
\hline
1&1&1  &   &    &   &   &   &  & 2\\
\hline
2&1 &4& 1 &  &  &  &  &  & 6\\
\hline
3&1&3&3&1   &  &  &  &  &8\\
\hline
4&7&32&12&32&7  &  &  & &90\\
\hline
5&19&75&50&50&75&19 &  &  & 288\\
\hline
6& 41& 216 & 27 & 272& 27 & 216 & 41 &  & 840\\
\hline
7&751&3577&1323&2989&2989&1323&3577&751& 17280\\
\hline
  \end{array}
  $$
  \caption{$Q_n\left( f^ {(1)} \right) = (z-x)/c_n\, \sum_{i=0}^ n A_i \, f^ {(1)} (\xi_i(x))$. \label{tabnc}}
  \end{table}

 \subsection{Newton-Rectangle Iterative Function}\label{subnewton}
 \medskip
 \noindent
For  $n=0$, the {\em left rectangle rule} uses  an unique node
(the left end of the integration interval). When  this rule is applied to the
integral $\int_x^ z f^ {(1)}( t)\, dt$ we obtain
$$
Q_0\left( f^ {(1)} \right)= (z-x) \, f^ {(1)} (x).
$$
In this case, the sum of the weights is $c_0=1$ and the function  $B_0(x)$ 
(defined by \eqref{q3}) has the form $B_0(x)=  f^ {(1)} (\xi_0(x))=  f^ {(1)} (x)$. 
If $z$ is a simple root of $f$, since $\xi_0(z)=z$, according to Proposition
\ref{conv}, the iterative method generated by
\begin{equation}\label{q15}
t_0(x)= x- c\, B_0^ {-1} (x) \, f(x)= x-\displaystyle{\frac{f(x)}{f'(x)}},
\end{equation}
converges to the fixed point $t_0(z)=z$, and the local convergence is superlinear. 
The mapping $t_0(x)$ is coincident with the Newton's iterative function.
 \subsection{Newton-Trapezoidal Iterative Function }\label{main}
  
   \medskip
 \noindent
 When $n=1$, the {\em  trapezoidal rule} uses as nodes both ends of 
 the integration interval. We can thus define the stepsize
   $h_1(x)$ satisfying
 $$
  h_1(x)= t_0(x)-x\quad \Rightarrow \quad h_1^ {(1)}(x)= t_0^ {(1)}(x)-1, \quad  h_1^ {(2)}(x)= t_0^ {(2)}(x),
  $$
where $t_0$ is defined by \eqref{q15}.

Applying the mentioned rule to $\int_{ x}^ z f^ {(1)}(t) dt$,  with  nodes $\xi_0(x)= x$ and $\xi_1(x)= x+h_1(x)$,  we obtain
  $$
  Q_1(f^ {(1)})= \displaystyle{ \frac{z-x}{c_1}   } \, B_1(x)=\displaystyle{ \frac{z- x}{2}   } \, \left[ f^ {(1)}( x)+ f^ {(1)}(x+ h_1(x)) \right].
  $$
  Therefore the iterative function has the form
  \begin{equation}\label{q17}
  \begin{array}{ll}
t_1(x)&= x- c_1 \, B_1^ {-1}(x) f(x)= x-   \displaystyle{ \frac{2\, f( x)}{f^ {(1)}(x)+ f^ {(1)}\left(x+h_1(x)\right)}   }=\\
 & =x-   \displaystyle{ \frac{ 2\,f( x)}{f^ {(1)}(x)+ f^ {(1)}\left(\displaystyle{x-\frac{f(x)}{f^ {(1)}(x)}}\right)}   } .
 \end{array}
\end{equation}
 The last formula can also be written as 
 \begin{equation}\label{q18}
 t_1(x)= x- \displaystyle{    \frac{f(x)}{  
\displaystyle{  \frac{f^ {(1)} (x)+ f^ {(1)} \left(x-f(x)/f^ {(1)}(x)\right)}{2}}  } 
       } .
\end{equation}
The equation \eqref{q18} means that the {\em step} of $t_1(x)$, that is $|H_1(x)|=|t_1(x)-x|$, is obtained from the average between the slopes of the tangents to the graphic $y=f(x)$ at the points $x$ and $\bar x$, where $\bar x=t_0(x)=x-f(x)/f^ {(1)}(x)$ (in other words, $\bar{x}$ is the image of $x$ by the iterative function $t_0$).

\medskip
\noindent
The function \eqref{q17} generates the iterative process
   \begin{equation}\label{q19}
    \begin{array}{l}
    h_k= t_0(x_k)-x_k=- f(x_k)/f^ {(1)}(x_k)\\
x_{k+1}=x_k-   \displaystyle{ \frac{2\, f(x_k)}{f^ {(1)}(x_k)+f^ {(1)}( x_k+h_k  )}  }, \qquad k=0,1,\ldots,
   \end{array}
  \end{equation}
which we  denominate Newton\--trapezoidal.

    \begin{proposition}\label{propotrap1}
   Assume that the real function $f$ is sufficiently regular in a neighborhood of a certain simpre root $z$ and 
    an initial approximation $x_0$ is chosen sufficiently close to $z$. Then,  the Newton\--trapezoidal method
   \eqref{q19} converges to $z$ and its convergence order is at least  $3$.
  \end{proposition}
    
  \medskip
  {\bf Proof.}
  The nodes of the quadrature rule are
 $\xi_0(x)=x$ and $\xi_1(x)= x+h_1(x)=t_0(x)$. Hence $\xi_0(z)=\xi_1(z)=z$.
 By  Proposition \ref{conv}, the point  $z$ is a superatractor fixed point of $f$, that is
 $t_1(z)=z$ and $t_1^ {(1)} (z)=0$.

\noindent
Since
$$
B_1(x)= f^ {(1)}(x)+ f^ {(1)} (x+h_1(x))=f^ {(1)}(x)+ f^ {(1)} (t_0(x)), 
$$
we have
$$
\begin{array}{ll}
B_1^ {(1)}(x)&=f^ {(2)}  (x)+ f^ {(2)} (x+ h_1(x)) \, (1+h_1^ {(1)} (x))     \\
&=f^ {(2)}  (x)+f^ {(2)}  (t_0(x))\, (1+t_0^ {(1)}(x)-1).
\end{array}
$$
Therefore, since $t_0^ {(1)}(z)=0$, we get
$$
B_1^ {(1)}(z)= 2\, f^ {(2)}  (z). 
$$
On the other hand, as $c_1=2$, from \eqref{q13} we obtain
$$
B_1^ {(1)}(x)\, H_1(x)+ B_1(x)\, H_1^ {(1)} (x) =- 2\, f^ {(1)} (x).
$$
Thus,
$$
B_1^ {(2)}(x)\, H_1(x)+  2\, B_1^ {(1)}(x)\, H_1^ {(1)} (x)+ B_1(x)\, H_1^ {(2)}(x) =- 2\, f^ {(2)} (x).
$$
Note that for $x=z$ we have $H_1(z)=0$ and $H_1^ {(1)} (z)= t_1^ {(1)}(z)-1 =-1$, yielding
$$
-2\, f^ {(2)}(z)+ 2 \, f^ {(1)} (z) \, H_1^ {(2)} (z) = -2 \, f^ {(2)} (z).
$$
Since  $z$ is a simple root of $f$, from the last equality we conclude that
$$
H_1^ {(2)} (z)=0 .
$$
Noting that $H_1 ^{(2)} (x) \equiv  t_1^ {(2)} (x)$, we finally obtain
$$
 t_1^ {(2)} (z)=0.
$$
Therefore the method
 \eqref{q19} converges locally to the root $z$ of  $f$ and its convergence order is at least 3.
$ \hfill \Box$

\medskip
   In the case of a multivariate function $f:D\subset \Rb^ d\mapsto \Rb^ d$, the iterative function of the
   Newton\--trapezoidal method can be written in the form
    \begin{equation}\label{qnt}
  \begin{array}{l}
  h_1(x)=\displaystyle{J_f^{-1} f(x)},\qquad  x \in \Rb^ d\\
t_1(x)=x-  2\, \left[ J_f(x)+  J_f(x+  h_1(x))  \right]^ {-1}\, f(x),
\end{array}
\end{equation}
where $J_f$ denotes the Jacobian matrix of the function $f$.

 \subsection{ Newton\--Simpson Iterative Function}\label{secsimp}

For  $n=2$,  applying the Simpson's rule to the integral $\int_{x}^ z f^ {(1)}(t)\, dt$ we obtain
\begin{equation}\label{q23}
\begin{array}{ll}
 Q_2(f^ {(1)})&=\displaystyle{ \frac{z-x}{c_2}   }\, B_2(x)=\\
 &=  \displaystyle{ \frac{z-x}{6}   } \, \left[ f^ {(1)}( x)+ 4\,  f^ {(1)}\left( x+ h_2(x)\right) +  f^ {(1)}\left( x+2\,  h_2(x)\right) \right].
 \end{array}
\end{equation}
In \eqref{q23} the step $h_2(x)$ is defined recursively by means of the Newton-trapezoidal iterative function $t_1(x)$ given by  \eqref{q17}, 
\begin{equation}\label{q24}
h_2(x)= \displaystyle{ \frac{t_1(x)-x}{n}} =    \displaystyle{ \frac{t_1(x)-x}{2}}  .
\end{equation}
The function $h_2$ and its derivatives satisfy:
\begin{equation}\label{q25}
\begin{array}{ll}
h_2^ {(1)}(x)&=\displaystyle{ \frac{t_1^ {(1)}(x)-1}{2} } \Longrightarrow h_2^ {(1)}(z)=- 1/2    \\
h_2^ {(2)}(x)&=\displaystyle{ \frac{t_1^ {(2)}(x)}{2} } \Longrightarrow h_2^ {(2)}(z)=0.    \\
\end{array}
\end{equation}
We designate the mapping
\begin{equation}\label{q26}
t_2(x)= x- c_2\, B_2^ {-1} (x) \, f(x)
\end{equation}
as Newton\--Simpson iterative function. The corresponding iterative method
can be described as 
    \begin{equation}\label{q27}
    \begin{array}{l}
    h_2(x_k)=\displaystyle{\frac{ t_1(x_k)-x_k}{2}}\\
x_{k+1}=x_k-   \displaystyle{ \frac{6\, f(x_k)}{f^ {(1)}(x_k)+4\, f^ {(1)}( x_k+h_2(x_h)  )+ f^ {(1)}\left( x_k+2\, h_2(x_k)  \right)}  }, \qquad k=0,1,\dots
   \end{array}
  \end{equation}
  
  \begin{proposition}\label{propotrap2}
  Let $f$ be a sufficiently regular real function on a neighborhhood of a simple root $z$. Taking an initial approximation $x_0$  sufficiently close to $z$, the Newton-Simpson method  \eqref{q27} converges to $z$ and its convergence order is not less than $4$.
  \end{proposition}
  {\bf Proof.}
  \noindent
By \eqref{q23}, we have $B_2(z)= c\, f^ {(1)} (z)\neq 0$ and
 \begin{equation}\label{q28}
 \begin{array}{ll}
 B_2^ {(1)} (x)&= f^ {(2)}(x)+ 4\, f^ {(2)} (x+h_2(x))\, (1+h_2^ {(1)} (x))+\\
 & + f^ {(2)} (x+ 2\, h_2(x))\, (1+2 \, h_2^ {(1)} (x)).
 \end{array}
  \end{equation}
  Since $h_2(z)=0$, using \eqref{q25} we obtain
  \begin{equation}\label{q29}
  \begin{array}{ll}
  B_2^ {(1)} (z)&= f^ {(2)} (z)+ 4\, f^ {(2)} (z)\, (1-1/2)+ f^ {(2)} (1-1)\\
  & =3\, f^ {(2)}(z) =\displaystyle{  \frac{c_2}{2}}\, f^ {(2)}(z).
  \end{array}
  \end{equation}
  Differentiating both sides of \eqref{q28} yields
  \begin{equation}\label{q30}
  \begin{array}{ll}
  B_2^ {(2)} (x)&= f^ {(3)}(x)+\\
  &+ 4\left[ f^ {(3)} (x+ h_2(x))\, \left(1+ h_2^ {(1)} (x)\right)^ 2 + f^ {(2)} (x+h_2(x)) \, h_2^ {(2)} (x) \right]+\\
  &+2\,  f^ {(3)} (x)\, h_2(x) \left(1+ 2\, h_2^ {(1)} (x)\right)^ 2+  f^ {(2)}\left( x+2 \, h_2(x)\right)\, h_2^ {(2)} (x).
  \end{array}
  \end{equation}
  Since $h_2(z)=0$, $h_2^ {(1)} (z)=-1/2$ and $h_2^ {(2)}(z)=0$, from \eqref{q30} we conclude that
  \begin{equation}\label{q30a}
  \begin{array}{ll}
    B_2^ {(2)} (z)&=  f^ {(3)}(z)+4\, f^ {(3)} (z)\times (1/4)\\
    &= 2\, f^ {(3)} (z) =\displaystyle{\frac{c_2}{3}} \, f^ {(3)} (z).
  \end{array}
  \end{equation}
  Concerning the function  $H_2(x)= t_2(x)-x$, from \eqref{q26} we get
  $$
  B_2(x)\, H_2(x)=-c_2\, f(x).
  $$
  Differentiating three times the last equality, we obtain
  \begin{equation}\label{q31}
  B_2^ {(1)}(x)\, H_2(x)+ B_2(x)\, H_2^ {(1)}(x) = -c_2\, f^ {(1)} (x),
  \end{equation}
    \begin{equation}\label{q32}
  B_2^ {(2)}(x)\, H_2(x)+ 2\, B_2^ {(1)}(x)\, H_2^ {(1)}(x) + B_2(x)\, H_2^ {(2)} (x)= -c_2\, f^ {(2)} (x),
  \end{equation}
  and
    \begin{equation}\label{q33}
    \begin{array}{ll}
  B_2^ {(3)}(x)\, H_2(x)+ 3\, B_2^ {(1)}(x)\, H_2^ {(2)}(x) +\\
  \hspace{2cm}+3\,  B_2^ {(2)}(x)\, H_2^ {(1)} (x)+ B_2(x)\, H_2^ {(3)} (x)= -c_2\, f^ {(3)} (x).
  \end{array}
  \end{equation}
  Since $H_2(z)=0$ and $B_2(z)= c_2 \, f^ {(1)} (z)\neq 0$, it follows from \eqref{q31} that $H_2^ {(1)} (z)=-1$,  that is, $t_2^ {(1)}(z)=0$, and therefore the corresponding iterative method has convergence order at least 2
  (as we know, from Proposition  \ref{conv}).
  
  \medskip
  \noindent
  From \eqref{q32} we obtain
  $$
  2\, B_2^ {(1)} (z)\, H_2^ {(1)} (z)+ B_2(z)\, H_2^ {(2)} (z) =- c_2\, f^ {(2)} (z),
  $$
  that is,
  $$
  -2 \, B_2^ {(1)} (z)+ B_2(z)\, H^ {(2)} (z)=- c_2\, f^ {(2)} (z).
  $$
  Taking  \eqref{q29} into consideration, we have
  $$
  -c_2\, f^ {(2)} (z)+ B_2(z) \, H_2^ {(2)} (z)= -c_2\, f^ {(2)} (z).
  $$
  Since $B_2(z)\neq 0$, we obtain $H_2^ {(2)}(z)=0$, that is, $t_2^ {(2)} (z)=0$, which means that the iterative method \eqref{q27} has convergence order not less than 3.
  
  \medskip
  \noindent
As $H_2^ {(1)} (z) = t_2^ {(1)}-1=-1$ and  $H_2^ {(2)} (z)= t_2^ {(2)} (z)=0$, from \eqref{q33} it follows that
  $$
  3\, B_2^ {(2)} (z)\, H_2^ {(1)} (z)+ B_2(z)\, H_2^ {(3)}(z)=-c_2\, f^ {(3)} (z),
  $$
  that is,
  $$
  -3\, B_2^ {(2)} (z)+ B_2(z) \, H_2^ {(3)} (z) = - c_2\, f^ {(3)} (z).
  $$
  Finally, from  \eqref{q30a} and \eqref{q11} we obtain
  $$
   -c_2\,  f^ {(3)} (z)+ c_2 \, f^ {(1)} (z)\, H_2^ {(3)} (z) = - c_2\, f^ {(3)} (z).
  $$
 Hence $H_2^ {(3)}(z)=0\Longleftrightarrow t_2^ {(3)} (z)=0$. Therefore we may conclude  that the method \eqref{q27} 
 has convergence order not less than 4.
 $\hfill \Box$

\medskip
  In the case of a multivariate function $f:D\subset \Rb^ d\mapsto \Rb^ d$, the iterative function of the
   Newton\--Simpson method can be written in the form
    \begin{equation}\label{q34}
  \begin{array}{l}
  h_2(x)=\displaystyle{\frac{t_1(x)-x}{2}},\qquad  x \in \Rb^ d\\
t_2(x)=x-  6\, \left[ J_f(x)+ 4\,  J_f(x+  h_2(x)) +  J_f(x+ 2\,  h_2(x))  \right]^ {-1}\, f(x),
\end{array}
\end{equation}
where  $t_1$ is the iterative function of the Newton-trapezoidal method in $\Rb^ d$, defined by \eqref{qnt}.

\medskip
{ \bf Remark.} {\it 
One can verify that if in  \eqref{q24} we replace $t_1$ by $t_0$ (that is, if we define $h_2(x)=(t_0( x)-x)/2$ instead of $h_2(x)=(t_1(x)-x)/2$), the resulting method has just third, and not fourth order of convergence. This confirms the  advantage of the recursive process we have introduced here to define the Newton-Simpson and the subsequent iterative functions.
}

\medskip

\noindent
For the sake of simplicity, in the rest of this paper we shall refer to the iterative methods corresponding to the functions $t_0$, $t_1$ and $t_2$ as Newton's, Trapezoidal and Simpson's methods, respectively.

\subsection{ Convergence order of Newton\--Cotes iterative functions }\label{ordem}

\noindent
Just as in the case of the Trapezoidal and Simpson's methods, in the general case, for $n\geq 3$, the convergence order of an iterative function $t_n$, based on a quadrature rule, depends only on $n$ and on the convergence order of $t_0$.
We will now prove some lemmas that will be used later to obtain the main result of this paper, the Theorem \ref{principal},
concerning the convergence order of the iterative functions $t_n$.

\medskip
 For $n\geq 1$, we assume that $z$ is a simple root of the real function $f$, which is differentiable, at least, $n$ times, in a neighborhood of $z$. Let $t_n$ be an iterative function based on a Newton-Cotes closed rule, defined by $t_n(x)=x- c_n\, B_n(x)^ {-1}\, f(x)$, with
$B_n(x)= \sum_{i=0}^n A_i\, f^ {(1)} (x+i\, h_n(x))$, where $h_n(x)= (t_{n-1}(x)-x)/n $ denotes the step of the quadrature rule and $c_n=\sum_{i=0}^ n A_i$ denotes the sum of its 
weights.
\begin{lemma}\label{basico}
 The step $h_n(x)=\displaystyle{ \frac{t_{n-1}(x)-x}{n} }$ of the iterative function $t_n(x)$, at  $x=z$,
  satisfies
\begin{equation}\label{new1}
\begin{array}{ll}
(i)& h_n(z)=0\\
(ii)& h_n^ {(1)} (z)=- 1/n\\
(iii)&   h_n^ {(j)} (z)=-\displaystyle{ \frac{t_{n-1}^ {(j)} (z)}{n} }\quad j=2,3,\ldots,n .
\end{array}
\end{equation}
\end{lemma}
{\bf Proof.}
From proposition
\ref{conv} we know that $t_{n-1}$ has superlinear convergence. Thus, $t_{n-1}(z)=z$ and $t_{n-1}^ {(1)} (z)=0$, 
from where the equalities (i) to (iii) follow, taking into consideration the definition of $h_n$.
$\hfill \Box$

\medskip
The $A_i$ coefficients in the function $B_n(x)$ 
satisfy certain equalities which are the subject of the following Lemma.

 \begin{lemma}\label{lema1}
Let  $A_0$, $A_1$, $\ldots, A_n$ be
the weights of a Newton-Cotes quadrature rule, with $n+1$ nodes. The following equalities hold:
\begin{equation}\label{ig}
\begin{array}{ll}
(i_1)& \displaystyle{ \frac{A_1+2\, A_2+ \cdots+n\, A_{n}}{n}    }= \displaystyle{ \frac{c_n}{2}  }\\
(i_2)& \displaystyle{ \frac{ A_1+ 2^ 2 \, A_2+ \cdots+ n^ 2\, A_n}{n^ 2}    }= \displaystyle{ \frac{c_n}{3}  }\\
\vdots&\hspace{3cm} \vdots\\
(i_n)& \displaystyle{ \frac{A_1+ 2^ n\, A_2+ \cdots+ n^ n\, A_n}{n^ n}    }= \displaystyle{ \frac{c_n}{n+1}  },
\end{array}
\end{equation}
where  $c_n=\sum_{i}^ n\, A_i$.
\end{lemma}
{ \bf Proof.}
We begin by considering the case where the integration interval is $[0,n]$.
In this case, using the method of undetermined coefficients, the weights $A_i$ are the solution of the following linear system,
\begin{equation} \label{sys}
\begin{array}{ll}
A_0+A_1+ \cdots+A_n   &= n\\
\qquad  A_1+2\, A_2+ 3\, A_3+\cdots+ n\, A_n &= n^ 2/2\\
\qquad  A_1+ 2^ 2\, A_2+ 3^2\, A_3+ \cdots+ n^ 2\, A_n&= n^ 3/3\\
& \vdots \\
\qquad  A_1+ 2^ n\, A_2+ 3^ n\, A_3+ \cdots+ n^ n \, A_n &=n^ {n+1}/(n+1).\\
\end{array}
\end{equation}
Since we have $c_n=n$, it is obvious that the last $n$ equations of the system \eqref{sys}
are equivalent to the system \eqref{ig} and therefore they have the same solution.
To deal with the case of an integration interval of arbitary length $c_n$, we take into account that
in this case the nodes are $x_i=i \,c_n/n$. If in \eqref{sys} we replace $i$ by  $x_i=i\, c_n/n$,
$i=0,...,n$, we again obtain a linear system which is equivalent to \eqref{ig}.
$\hfill \Box$

\medskip
The equalities \eqref{ig} are used in the proof of the following Lemma.

\begin{lemma}\label{lema2}
Given  $n\geq 1$, let $z$ be a simple root of a real function $f$, differentiable up to order  $n$ at least, in a neighborhood of
$z$. Let $t_0$ be the  Newton's  function and  $t_n$ be the  Newton\--Cotes closed iterative function with $n+1$ nodes, $t_n(x)=x- c_n\, B_n(x)^ {-1}\, f(x)$, with
$B_n(x)= \sum_{i=0}^n A_i\, f^ {(1)} (x+i\, h_n(x))$,  where $h_n(x)=(t_{n-1}(x)-x)/n$ and $\sum_{i=0}^ n\, A_i=c_n$.

\noindent
If $t_{n-1}$ has convergence order at least $n$, that is,
$$
t_{n-1}(z)=z\quad \mbox{and} \quad  t_{n-1}^ {(1)} (z)=t_{n-1}^ {(2)}(z)= \ldots = t_{n-1}^ {(n-1)}=0,
$$
then  $B_n$ and its  $n$  first derivatives satisfy, at $x=z$, the following equalities:
\begin{equation}\label{ig2}
B_n^ {(j)} (z)= \displaystyle{ \frac{c_n}{j+1}  }\, f^ {(j+1)} (z),\quad\mbox{para}\quad  j=0,1,2,\ldots,n.
\end{equation}
\end{lemma}
{\bf Proof.}
The proof will  use induction on $j$.
We first note that 
$$B_n(x)=B_n^{(0)}(x)=\sum_{i=0}^ n A_i\, f^ {(1)}(x+i h_n(x)).$$
 Since $h_n(z)=0$ (see \eqref{new1}), we conclude that 
 $$B_n(z)=\sum_{i=0}^n A_i\, f^ {(1)} (z)=c_n\, f^ {(1)} (z),$$
 which means that \eqref{ig2} is true for $j=0$.
\medskip
\noindent
We shall now deal with the case $j\geq 1$.
First we will show that
\begin{equation}\label{bnj}
B_n^{(j)}(x)=\sum_{i=0}^n A_i\left[  f^{(j+1)} (x + i \,h_n(x) ) (1 + i\, h_n^{(1)}(x))^j + \cdots \right],
\qquad j=1,2, \dots
\end{equation}
where the omitted terms on the right-hand side of \eqref{bnj} contain derivatives of $h_n$ of order greater
than one.
For $j=1$, differentiation of $B_n$ gives
\begin{equation}\label{bn1}
B_n^{(1)}(x)=\sum_{i=0}^n A_i\left[  f^{(2)} (x + i h_n(x) ) (1 + i h_n^{(1)}(x)) \right],
\end{equation}
which is in agreement with \eqref{bnj}.
Suppose now, as induction hypothesis, that \eqref{bnj} is true for $j=1,2,..,k$. Let us show that it also holds
for $j=k+1$. We have
\begin{equation}\label{bnk}
B_n^{(k+1)}(x)=(B_n^{(k)}(x))^{(1)}=\sum_{i=0}^n A_i\left[  f^{(k+1)} (x + i\, h_n(x) ) (1 + i \,h_n^{(1)}(x))+ \cdots \right]^{(1)}.
\end{equation}
Hence
\begin{equation}\label{bnk1}
\begin{array}{c}
B_n^{(k+1)}(x)=\sum_{i=0}^n A_i\left[  f^{(k+2)} (x + i \,h_n^{(1)}(x) ) (1 + i \,h_n^{(1)}(x))^{k+1}+ \right.\\
 \left. + f^{(k+1)} (x + i \,h_n(x) ) k (1 + i\, h_n^{(1)}(x))^{k-1} h_n^{(2)}+\cdots \right].
\end{array}
\end{equation}
Since the term which contains $f^{(k+1)}$, on the right-hand side of \eqref{bnk1}, includes the second derivative of $h_n$, this term can be omitted, yielding
\begin{equation}\label{bnk2}
B_n^{(k+1)}(x)=\sum_{i=0}^n A_i\left[  f^{(k+2)} (x + i h_n^{(1)}(x) ) (1 + i h_n^{(1)}(x))^{k+1}+ \dots \right].
\end{equation}
Therefore, by mathematical induction, we conclude from \eqref{bnk2} that \eqref{bnj} holds for any natural $j$.

If in \eqref{bnj} we take the limit as $x \rightarrow z$, taking into account that $h_n(z)=0$, $h_n^{(1)}(z)=-1/n$,
and $h_n^{(j)}(z)=0$, for $j\ge 2$, we obtain
\begin{equation}\label{bnk2a}
B_n^{(j)}(z)=\sum_{i=0}^n A_i f^{(j+1)} (z) (1 -i/n)^{j}.
\end{equation}
 To complete the proof of Lemma 3.3 it remains to show that
 \begin{equation}\label{bnk3}
 \sum_{i=0}^n A_i \left( 1-i/n  \right)^j=\frac{c_n}{j+1}.
\end{equation}
Rewriting \eqref{ig} in the form
 \begin{equation}\label{bnk4}
 \sum_{i=0}^n A_i \left( i/n  \right)^j=\frac{c_n}{j+1}, \quad j=1,\dots,n
\end{equation}
and taking into account the symmetry $A_{n-i}=A_i$, for $j=0,1,\ldots,n$, we obtain
 \begin{equation}\label{bnk5}
 \sum_{i=0}^n A_i \left( \frac{n-i}{n} \right)^j=\frac{c_n}{j+1}, \quad j=1,\ldots,n,
\end{equation}
which is equivalent to \eqref{bnk3}.
$\hfill \Box$

\bigbreak
The main result follows.
\begin{thm}\label{principal}
Let $z$ be a simple zero of the real function $f$ continuously differentiable, up to order $n+2 $, in a neighborhood of $z$.
For $n\geq 1$,  the Newton-Cotes iterative function $t_n$, defined recursively from $t_0$, has local order of convergence  at least $n+2$.
\end{thm}

{\bf Proof.} The proof is by induction on $n$.
From Proposition 2.1, we know that $t_0$ has at least quadratic convergence and  by
Propositions 3.1 and 3.2, the iterative functions $t_1$ and $t_2$ have order of convergence at least 3 and 4, respectively. Let us suppose that for a certain $n\geq 2$, $t_n$ has convergence order at least $n+2$, that is, 
$$
t_n(z)=z,\quad t_n^ {(1)}(z)=t_n^ {(2)}(z)=\ldots=t_n^ {(n+1)}(z)=0.
$$
We need to prove that
$$
t_{n+1}(z)=z,\quad t_{n+1}^ {(j)}(z)=0, \quad\mbox{with}\quad 1\leq j\leq n+2.
$$
Let $t_{n+1}(x)= x- c_{n+1}\, B_{n+1}^ {-1}(x)\, f(x)$ and $H_{n+1}(x)=t_{n+1}(x)-x$. Then
\begin{equation}\label{AA}
B_{n+1}(x)\, H_{n+1}(x)= - c_{n+1}\, f(x).
\end{equation}
Note that $B_{n+1}(z)\neq 0$, since $B_{n+1}(z)= c_{n+1}\, f^ {(1)} (z)$ (see \eqref{ig2} ) and $z$ is a simple root. Thus, we conclude from  \eqref{AA} that
\begin{equation}\label{BB}
H_{n+1}(z)=0, \quad \mbox{i.e.}\quad t_{n+1}(z)=z.
\end{equation}
Moreover
\begin{equation}\label{CC}
\begin{array}{ll}
(i)& H_{n+1}^ {(1)} (z)= t_{n+1}^ {(1)} (z)-1\\
(ii)& H_{n+1}^ {(j)} (z)= t_{n+1}^ {(j)} (z), \quad j\geq 2.\\
\end{array}
\end{equation}
Let us differentiate both sides of \eqref{AA}, applying the Leibniz rule to the left-hand side.
 For $1\leq k\leq n+1 $ we obtain:
\begin{equation}\label{hello}
\begin{array}{ll}
(i_1)& B_{n+1}^ {(1)}\, H_{n+1} (x)+ B_{n+1}(x)\, H_{n+1}^ {(1)} (x) = -c_{n+1}\, f^ {(1)} (x)\\
(i_2)& B_{n+1}^ {(2)}\, H_{n+1} (x)+\binom{2}{1}\, B_{n+1}^ {(1)} (x)\, H_{n+1}^ {(1)} (x)+ B_{n+1} (x)\, H_{n+1}^ {(2)} (x)= -c_{n+1}\, f^ {(2)} (x)\\
(i_3)& B_{n+1}^ {(3)}\, H_{n+1} (x)+\binom{3}{1}\, B_{n+1}^ {(2)} (x)\, H_{n+1}^ {(1)} (x)+\\
& \qquad \hspace{2cm} + \binom{3}{2}\, B_{n+1}^ {(1)} (x)\, H_{n+1}^ {(2)} (x)+ B_{n+1} (x)\, H_{n+1}^ {(3)} (x)= -c_{n+1}\, f^ {(3)} (x)\\
&\hspace{3cm} \vdots\\
(i_{n+1})& \sum_{i=0}^ p \binom{p}{i}\, B_{n+1}^ {(n+1-i)} (x)\, H_{n+1}^ {(i)} (x)=-c_{n+1}\, f^ {(n+1)} (x).
\end{array}
\end{equation}
Taking  \eqref{BB} into account, in the equalities \eqref{hello}, from $(i_1)$ to $(i_{n+1})$, all the terms  containing $H_{n+1}(z)$ vannish when $x$ is replaced by $z$. Moreover, we know from \eqref{ig2}  that $B_{n+1}(z)= c_{n+1}\, f^ {(1)} (z)$.
Therefore, we can rewrite \eqref{hello} $(i_1)$ as 
$$
c_{n+1}\, f^ {(1)} (z)\, H_{n+1}^ {(1)}(z)= -c_{n+1}\, f^ {(1)} (z).
$$ 
Thus,
\begin{equation}\label{DD}
H_{n+1}^ {(1)} (z)= -1.
\end{equation}
From  \eqref{CC} it follows that
$$
t_{n+1}^ {(1)} (z)=0.
$$
Taking  \eqref{DD} into consideration, the equality  \eqref{hello} $(i_2)$  can be rewritten as
\begin{equation}\label{EE1}
-2 \, B_{n+1}^ {(1)} (z)+ B_{n+1}\, H_{n+1}^ {(2)} (z) = - c_{n+1} \, f^ {(2)} (z).
\end{equation}
From \eqref{ig2}  and \eqref{EE1} we conclude that 
$$
-2 \, B_{n+1}^ {(1)} (z)+ c_{n+1}\, f^ {(1)} (z) \, H_{n+1}^ {(2)} (z) = -c_{n+1}\, f^ {(2)}(z),
$$
and, by \eqref{ig2},
$$
-2 \,\displaystyle{  \frac{c_{n+1}}{2}   } f^ {(2)} (z)+ c_{n+1}\, f^ {(1)} (z) \, H_{n+1}^ {(2)} (z) = -c_{n+1}\, f^ {(2)}(z).
$$
Hence we get
\begin{equation}\label{FF}
H_{n+1}^ {(2)} (z)=0 \, \Longrightarrow t_{n+1}^ {(2)} (z)=0.
\end{equation}

To show that $t_{n+1}^ {(k)} (z)=0 $, for $k=3,...,n+2$ we use similar arguments.
First we rewrite \eqref{hello} $(i_{k})$, with $x=z$, taking into account that $H_{n+1}^ {(1)} (z)=-1$,
$ H_{n+1}^ {(2)} (z)=\dots= H_{n+1}^ {(k-1)} (z)=0$. We obtain
 $$
-k \, B_{n+1}^ {(k-1)} (z)+B_{n+1} (z) \, H_{n+1}^ {(k)} (z) = -c_{n+1}\, f^ {(k)}(z).
$$
Then, from \eqref{ig2}  it follows that 
$$
-k \,\displaystyle{  \frac{c_{n+1}}{k}   } f^ {(k)} (z)+ c_{n+1}\, f^ {(1)} (z) \, H_{n+1}^ {(k)} (z) = 
-c_{n+1}\, f^ {(k)} (z).
$$
Finally, from the last equation we conclude that
\begin{equation}\label{HH}
H_{n+1}^ {(k)} (z)=0  \, \Longrightarrow t_{n+1}^ {(k)} (z)=0.
\end{equation}
Since the last equality holds for $ 2\leq k\leq n+2$, the iterating function $t_{n+1}$ has order of convergence at least $n+3$.
This concludes the proof by induction.
$
\hfill \Box
$
\begin{remark}
If $z$ is a multiple root of $f$, the Theorem \ref{principal} cannot be directly applied. 
However, as referred in paragraph 2.1,  we can deal with this case by transforming the  
original equation  $f(x)=0$ into the equivalent equation $F(x)=0$,  such that $z$ is a simple root of $F$. 
Then we can construct iterative functions $t_n$ for $F$ and the Theorem \ref{principal} is applicable to them.
Example  4.3 in the next section illustrates this case.
\end{remark}

\begin{remark}
Note that in some cases the convergence order of $t_{n}$ can be higher than $n+2$
(see Examples 4.1 and 4.3 where, for even $n$, the convergence order of $t_n$ is $n+3$).
In such cases two consecutive iterative functions may have the same convergence order
(in the cited examples, the iterative functions $t_{2i}$ and $t_{2\,i+1}$ have convergence
order $2i+3$, $i=0,1,2,\ldots$).

\end{remark}
\section{Examples}\label{exemplos}

\noindent
Iterative processes with high order convergence may be particularly useful in cases where the choice of
{\em sufficiently close} initial approximations $x_0$ for the Newton's method is a difficult task.

As an illustration, consider  $f(x)= tanh(x-1)$ (Example \ref{exemplo1}).
The graphic of this function has the shape of a long flat $S$; in this case, the use of the Newton's method requires that the initial approximations $x_0$ be very close to the root (otherwise, the derivative of $f$ becomes very close to zero and the
Newton's method does not work).
It is worth to note that if we use the Newton-trapezoidal and the Newton-Simpson iterative functions,
we obtain convergence order $3$ and $4$, repectively, even with such initial approximations $x_0$ for which
the Newton's method  can not be applied.

\begin{figure}[h]
\begin{center} 
 \includegraphics[totalheight=5.2cm]{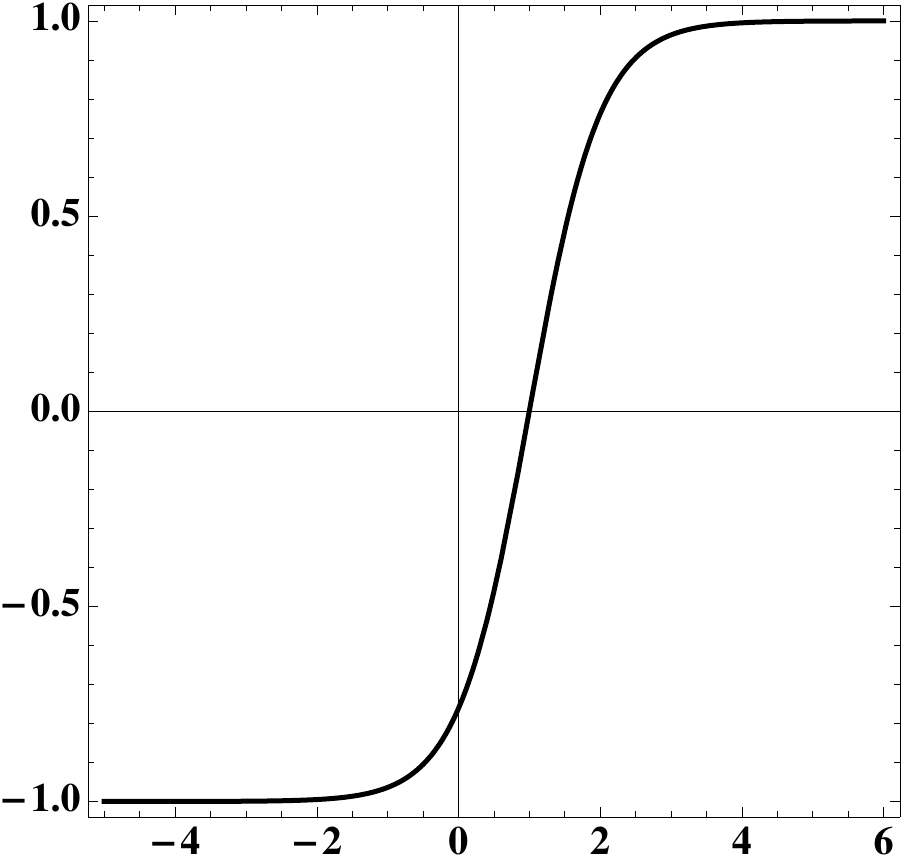}
\caption{\label{figtanh} $f(x)=tanh(x-1),\quad -5\leq x \leq 6$.}
\end{center}
\end{figure}

\begin{exemplo}\label{exemplo1}
 Fig.~\ref{figtanh} shows the graphic of the function 
$$
f(x)= tanh(x-1),
$$
in the interval $I=[-5,6]$.
  \end{exemplo}
The considered function $f$, which is infinitely differentiable,  has the unique root $z=1$ in this interval.
However, since the graphic has the shape of a long $S$, if we choose $x\in A=[-5,-0.9]$ or $x\in B=[2.9,6]$, we have
$|f(x)/f^ {(1)}(x)|>11$, which means that te Newton's iterative function $t_0$ has a long \lq\lq{step}\rq\rq\  $|t_0(x)-x|$.
Therefore, if the initial approximation satisfies $x_0\in A$ or $x_0\in B$ the subsequent iterates of the Newton's method get out of the interval $I$.

\medskip
\noindent
In this particular case, it is easy to verify that the Newton's method exceptionally  has cubic convergence,
since $t_0(1)=1$, $t_0^ {(1)}(1)=t_0^ {(2)}(1)=0$ and $t_0^ {(3)}(1)\neq 0$.
It may be verified that   Newton-trapezoidal iterative function $t_1$ has also convergence order $3$
and the  Newton-Simpson function $t_2$ has convergence order $5$.
\medskip
\noindent
In particular the Newton\--Simpson iterative function works when the initial approximation $x_0$ is chosen in a larger
interval (compared with the Newton's method), which follows from the fact that $t_2$ has a higher convergence order.

\noindent
In Table \ref{tab1} the values of the first derivatives of each iterative function, at $z=1$,
are displayed for comparison.
   \begin{table}
$$
\begin{array}{| c | c || c| c | c | c| c| c | }
\hline
t_i&t_i^ {(0)} & t_i^ {(1)} &t_i^ {(2)} & t_i^ {(3)} & t_i^ {(4)} &t_i^ {(5)} & \mbox{ordem} \\
\hline
t_0& 1 & 0 & 0 & -4 & 0 & -16 & 3\\
\hline
t_1& 1 & 0 & 0 & -1 & 0 & 14 & 3\\
\hline
 t_2& 1 & 0 & 0 & 0 & 0 & 82/3 & 5\\
\hline
  \end{array}
  $$
  \caption{Comparison of the convergence order of the Newton's method ($t_0$), Newton-trapezoidal method($t_1$)
   and Newton-Simpson ($t_2$), for Example 4.1 . \label{tab1}}
  \end{table} 
  
\begin{figure}[h]
\begin{center} 
 \includegraphics[totalheight=6.2cm]{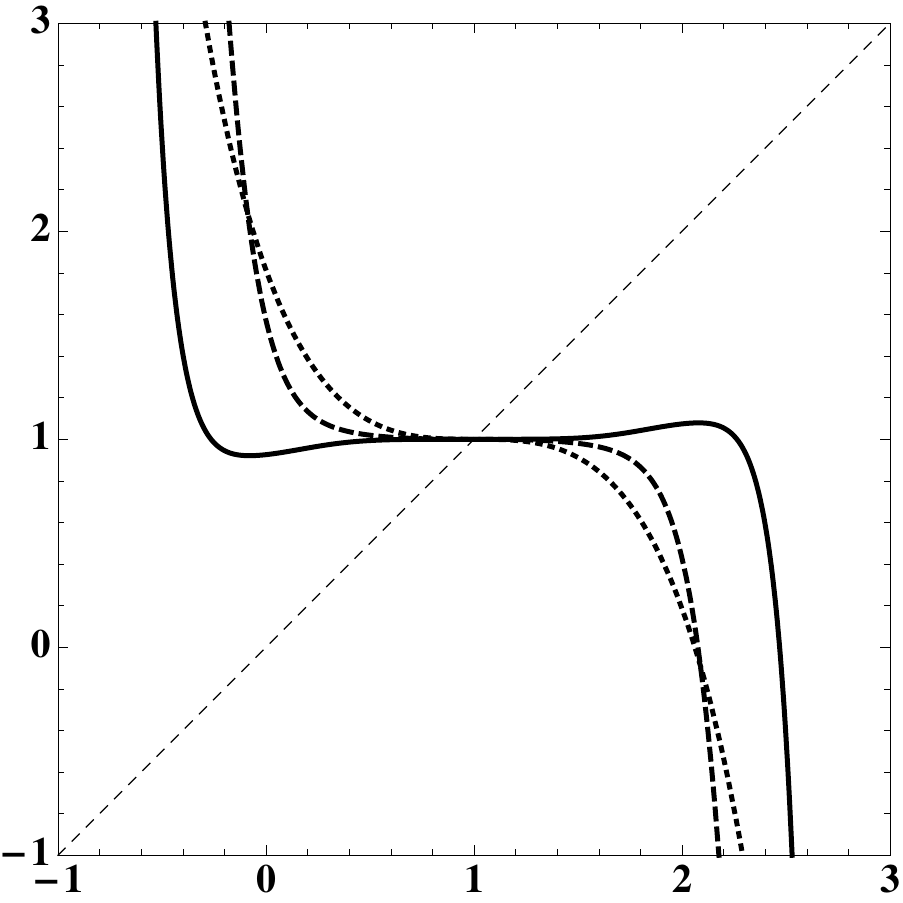}
\caption{\label{fignt1} Example 4.1: iterative functions $t_0$ (bold points), $ t_1$ (dashed line), and $t_2$ (full line).}
\end{center}
\end{figure}

  \medskip
  \noindent
 The graphics of the functions  $y=x$, $y=t_0(x)$ (Newton), $y=t_1(x)$ (trapezoidal ) and $y=t_2(x)$ (Simpson)
are displayed in Fig.  \ref{fignt1}.

\noindent
Though the trapezoidal method, in this case, has the same convergence order as the Newton's method, note that in the neighborhood of $z=1$ the graphic of $t_1$ is {\em flatter} than the one of $t_0$ \footnote{ We use here the term flat 
with a geometric intuitive sense, meaning {\em almost constant}. A more precise definition of this term will be given elsewhere.}.
\noindent
Analogously, since the graphic of $t_2$ (full line) starts to grow fast later than the other iterative functions, we conclude  that 
when using the corresponding method the initial approximation $x_0$ may be at a greater distance from $z=1$ than in the case of the Newton´s method, and that a small number of iterations of $t_3$ can produce a more accurate approximation of $z$, compared
with the result obtained with $t_0$. Using the same terminology as in \cite{rheinboldt}, p. 43, the Simpson's method has 
a larger { \em atraction basin} than the one of the Newton's method. The advantage of using methods whose 
{ \em atraction basin} is greater, specially in the context of numerical optimization without constraints,  will be discussed in detail in a future work.

\begin{figure}[h]
\begin{center} 
 \includegraphics[totalheight=4.2cm]{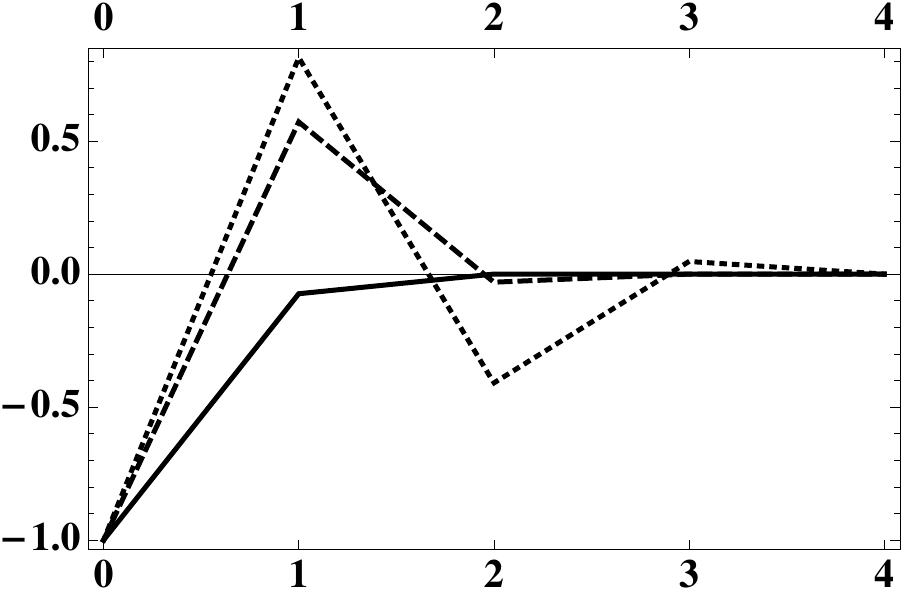}
\caption{\label{fig2}  $x_0=2.0$ (Example 4.1), error of $4$ iterates of  $t_0$ (points) , $ t_1$ (dashed line), and $t_2$ (full line).}
\end{center}
\end{figure}

\medskip
\noindent
Starting with $x_0=2.0$, four iterations have been computed for Newton's, trapezoidal and Simpson's method.  In Fig.~\ref{fig2} the error of the successive iterates is compared for the three methods (the computations were carried out using the {\sl Mathematica } \cite{wolfram1}
system with machine precision, that is, approximately 16 decimal digits). Note that the advantage of the Simpson's method,  in terms of accuracy,  compared with the other two methods, is visible from the first iteration onwards.

\begin{figure}[h]
\begin{center} 
 \includegraphics[totalheight=4.2cm]{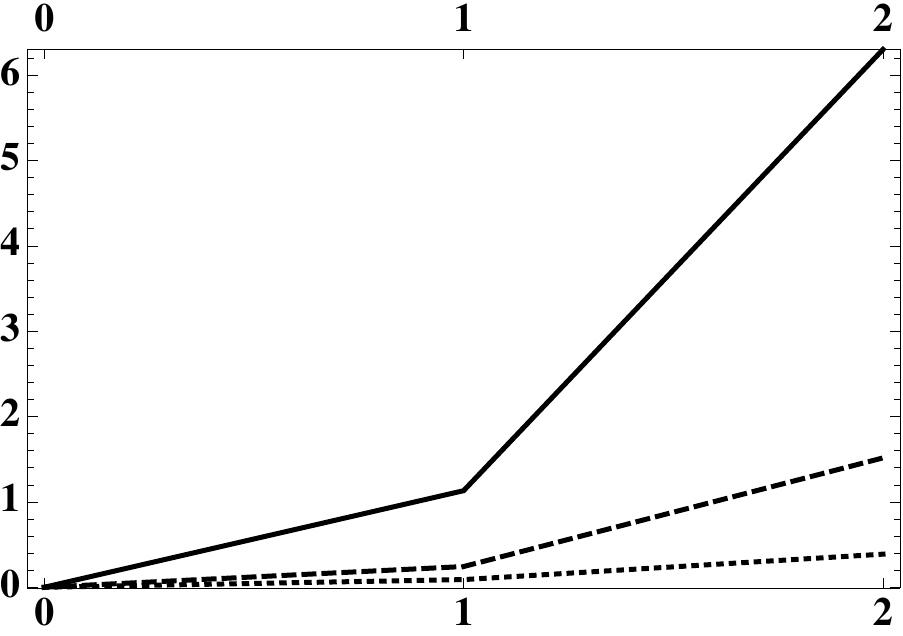}
\caption{\label{fig3} Example 4.1: $x_0=2.0$, number of significant digits for $2$ iterations of $t_0$ (points) , $ t_1$ (dashed), and $t_2$ (full).}
\end{center}
\end{figure}

\medskip
\noindent
Fig.~\ref{fig3} illustrates the improvement of accuracy which is obtained when an iterative function with convergence order $5$
is applied (which is the case of the Newton-Simpson process in this example). In this figure we show the  number $s$
of significant digits (that is, $s=- \log_{10} (| z-x_k|)$), corresponding to the two first iterates of the three mentioned
 methods, with $x_0=2.0$.

\noindent
An  even more impressive improvement of accuracy can be observed if iterative Newton-Cotes functions of higher order are used.
Consider, for example,  $t_4$, whose convergence order is $p+n=3+4=7$. In this case, the iterative function $t_4$ is given by:
$$
\begin{array}{ll}
h_4(x)&= (t_3(x)-x)/4\\
B_4(x)&= 7\, f^ {(1)} (x)+ 32\, f^ {(1)} \left(x+h_4(x)\right)+ 12\, f^ {(1)}\left( x+2\, h_4(x) \right)+\\
&+  32\, f^ {(1)} \left( x+ 3\, h_2(x) \right)+ 7\,  f^ {(1)} \left(x+ 4\, h_4(x)\right)\\
t_4(x)&= x- 90\, B_4^ {-1} (x) \, f(x).
\end{array}
$$

\begin{figure}[h]
\begin{center} 
 \includegraphics[totalheight=4.2cm]{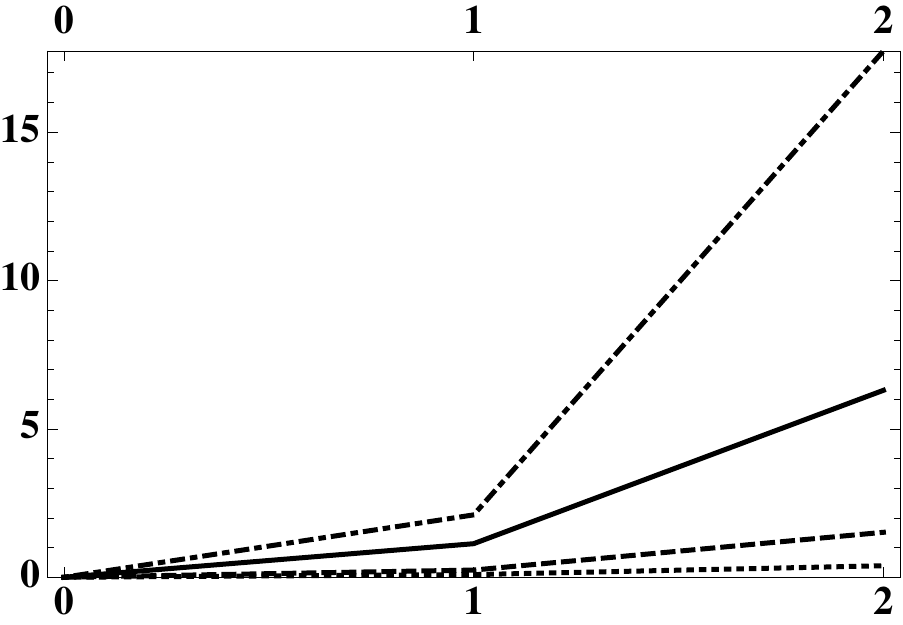}
\caption{\label{fig4A} Example  4.1:  number of significant digits for the two first iterations, with $x_0=2.0$, in the case of
 $t_0$ (points), $ t_1$ (dashed line),  $t_2$ (full line) and $t_4$ (dashed\--point).}
\end{center}
\end{figure}

\noindent
The first nonzero derivative at $z=1$ is $t_4^{(7)}(1) \approx -4.9$. In this case, the second iterate 
of the corresponding iterative method has already more than $17$ significant digits (see Fig.  \ref{fig4A}).

\medskip
\noindent
It is interesting to observe the numerical effect of {\em a single iteration} of each Newton-Cotes method, from $t_0$ to $t_7$.
With this purpose, we have chosen the initial approximation $x_0=1.1$, which we consider { \em sufficiently close}
to $z=1$, in the sense that with $x_0 = 1.1$  all the mentioned methods converge to $z$. The improvement of accuracy
after one iterate is shown in Table \ref{tab1nn}, where the number $s$ of significant digits is displayed, as well as the 
theoretical convergence order $q$ of each method. Note that $q=3+n-1$, when $n$ is odd, and $q=3+n$, when $n$ is even.
 \begin{table}
$$
\begin{array}{| c || c | c| c | c| c | c | c | c|}
\hline
t_i&t_0& t_1& t_2 & t_3 & t_4 & t_5 & t_6 & t_7\\
\hline
s&3.2&3.8&5.6&7.8&10.2&11.1&13.5&14.5\\
\hline
q&3&3&5&5&7&7&9&9\\
\hline
  \end{array}
  $$
  \caption{Results after 1 iteration with $x_0=1.1$ (Example 4.1). The parameter $s$ represents the approximate number of significant digits, and
$q$ is the convergence order.}   \label{tab1nn}
  \end{table}
\begin{table}
$$
\begin{array}{| c || c | c| c | c| c | c | }
\hline
t_{ij}&t_{21}& t_{32}& t_{43} & t_{54} & t_{65}& t_{76}\\
\hline
s&19.5&30.8&57.5&75.2&104.7&127.3\\
\hline
q&15&25&35&49&63&81\\
\hline
  \end{array}
  $$
  \caption{ Results after 1 iteration  with $x_0=1.1$ (Example 4.1). $s$ is the approximate number of significant digits.  
  $q $ is the convergence order.\label{tab1nnA}}
  \end{table}
  \begin{table}
$$
\begin{array}{| c || c | c| c | c| c | c | }
\hline
t_{ij}&t_{12}& t_{23}& t_{34} & t_{45} & t_{56}& t_{67}\\
\hline
s&17.7&39.5&53.4&80.9&98.8&135.4\\
\hline
q&15&25&35&49&63&81\\
\hline
  \end{array}
  $$
  \caption{Results after 1 iteration  with $x_0=1.1$ (Example 4.1). $s$ is the approximate number of significant digits and 
  $q $ is the convergence order.\label{tab1nnB}}
  \end{table}

 \medskip
  \noindent
  We have also applied some other methods, which result from the composition of two iterative functions, called
  composed methods.
  Let us denote 
  $$t_{ij}(x)= t_i(t_j(x)).$$ 
  Note that  we have $q(t_{ij}) =q(t_{i}) q(t_{j})$ (the convergence order
  of a composed method is the product of the convergence orders of the two components).
In Tables   \ref{tab1nnA} and \ref{tab1nnB} we compare the accuracy of a certain number of composed methods.
 
 \noindent
 Though the methods $t_{ij}$ and $t_{ji}$ have the same convergence order, we observe that the number $s$
 of significant digits, after one iteration,  is different in each case.
  
  \medskip
  \noindent
  When writing the code for the iterative Newton-Cotes functions in {\sl Mathe\-ma\-tica} we have used dynamical programming.
  Therefore, for example, once the value $x_1=t_6(x_0)$ is computed, the value $t_{76}(x_0)=t_7(x_1)$ can be obtained with a small additional effort. However this yields a very significant improvement of accuracy: from $13.5$ digits in the case
  of $t_6$ (see Table \ref{tab1nn}) to $127.3$ digits in the case of  $t_{76} $ (see Table \ref{tab1nnA}).

\begin{exemplo}\label{exemploA}

An extremal case of \lq\lq{bad behaviour}\rq\rq \ of the Newton's method occurs when, for a certain initial approximation $x_0\neq z$,  the sucessive iterates $x_{k+1}=t_0(x_k)$ are further and further apart from $z$. This happens when $z$ is a  {\em repelling} fixed point for $t_0$. For example, in the case of the (unique) fixed point $z=0$ of the function  \cite{benisrael1}, \cite{ujevic},
$$
f(x)=x^ {1/3}
$$
the derivative  $f^ {(1)}$ is not defined at $z=0$ and
$$
\lim_{x\rightarrow 0^ {+}} f^ {(1)}(x)= +\infty.
$$
Since $t_0(x)=x-f(x)/f^ {(1)}(x)=-2\, x\Longrightarrow  t_0^ {(1)}(0)=-2<-1$, and so
 $z=0$ is a repelling fixed point for $t_0$.
 \end{exemplo}

\noindent
If we consider the application of the Newton-Cotes iterative functions $t_i(x)$, with $i>0$, we come
to a similar conclusion, that is, $z=0$ is a repelling fixed point for all these functions.

\medskip
\noindent
In this case, we may apply the  procedure  suggested in paragraph \ref{secmult} for the case of multiple roots. That is, we may consider the equivalent equation  $F(x)=0$, with $F(x)=t_0(x)-x=-3\, x$. If we do so,
the corresponding iterative functions $t_i(x)$, starting with $i=0$ (Newton's method) are such that 
 $ t_i(x)=0$, $ \forall x \in \Rb$, $i \geq 0$. This means that we obtain the exact solution with the
 first iteration, for any initial approximation. 
 
 \noindent
 In conclusion, with this simple transformation of the equation,
 from an extremely difficult problem we obtain  an extremely easy one.

\begin{exemplo}\label{exemplo2}
The real function
$$
f(x)=\sin(x)-x,
$$
 has the unique real root $z=0$. However, this is a multiple root since $f(0)=0$ and $f^ {(1)}(0)=0$.
 Therefore, the Newton's method has local convergence order $p=1$.
\end{exemplo} 
 
 \noindent
 As can be seen from Table \ref{tabnova1},
 the performance of the Newton-Cotes iterative methods $t_1$ to $t_7$
in this case  is similar to the one of the Newton's method, that is, they don't offer any significant 
 advantage compared to $t_0$.

\begin{figure}[h]
\begin{center} 
 \includegraphics[totalheight=4.3cm]{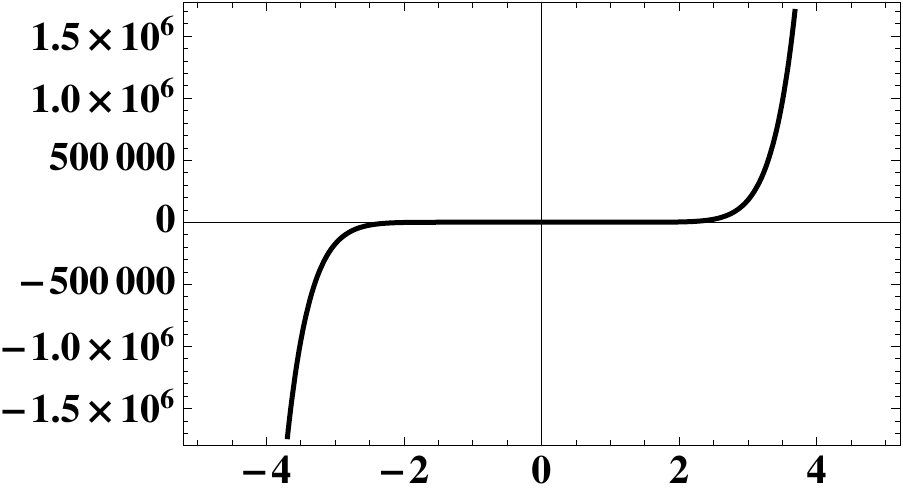}
\caption{\label{figpol1} $f(x)=x^ {11}+ 4\, x^ 2-10$.}
\end{center}
\end{figure}

  \begin{table}
$$
\begin{array}{| c || c | c| c | c| c | c | c | c|}
\hline
t_i&t_0& t_1& t_2 & t_3 & t_4 & t_5 & t_6 & t_7\\
\hline
s&1.18&1.27&1.28&1.35&1.41&1.45&1.49&1.52\\
\hline
  \end{array}
  $$
  \caption{$f(x)=\sin(x)-x$. Results of the first iteration with $x_0=0.1$.  \label{tabnova1}}
  \end{table}
 As suggested in Section \ref{secmult},  let us replace $f$ by the function
   $$
  F(x)=t_0(x)-x=-f(x)/f'(x)= \displaystyle{\frac{x-\sin(x)}{1-\cos(x)}},
  $$
 which may be extended to $x=0$ , with $F(0)=0$. Since $z=0$ is a {\em simple root} of $F$, when the Newton's method is applied to this function it has quadratic convergence. Actually, we have $F(0)=0$ and 
    $$
  F^ {(1)}(0)=\lim_{\delta\rightarrow 0} \displaystyle{ \frac{F(\delta)- F(0)}{\delta}  }=  \displaystyle{  \frac{1}{3}}\neq 0.
  $$
  Thus, if we apply the Newton-Cotes iterative functions $t_0$ to $t_7$ to the equation $F(x)=0$, we obtain the results displayed in Table \ref{tabnova2}.
  \begin{table}
$$
\begin{array}{| c || c | c| c | c| c | c | c | c|}
\hline
t_i&t_0& t_1& t_2 & t_3 & t_4 & t_5 & t_6 & t_7\\
\hline
s&4.2&4.8&7.6&9.6&13.1&14.2&17.7&18.7\\
\hline
q&3&3&5&5&7&7&9&9\\
\hline
  \end{array}
  $$
  \caption{$F(x)= \displaystyle{\frac{x-\sin(x)}{1-\cos(x)}}$ (Example 4.3). Results of the first iteration with $x_0=0.1$. The theoretical convergence order is   $q=2+n-1$ (if $n$ is odd) and $q=2+n$ (if $n$ is even).\label{tabnova2}}
  \end{table}
\begin{exemplo}\label{exemplo1nt}
Let
$$
f(x)=x^ {11}+ 4\, x^ 2-10.
$$
\end{exemplo}
Since the term $x^{11}$ is strongly dominant for the polynomial function $f$, the graphic of this function suggests 
the existence of a multipple root at $z=0$ (see Fig. \ref{figpol1}). However this isn't the case;
the considered polynomial has a single root $z\simeq 1.1$, which is the unique real root, and the Newton's method 
has local convergence order $p=2$ when applied to this function.

\begin{figure}[h]
\begin{center} 
 \includegraphics[totalheight=6.5cm]{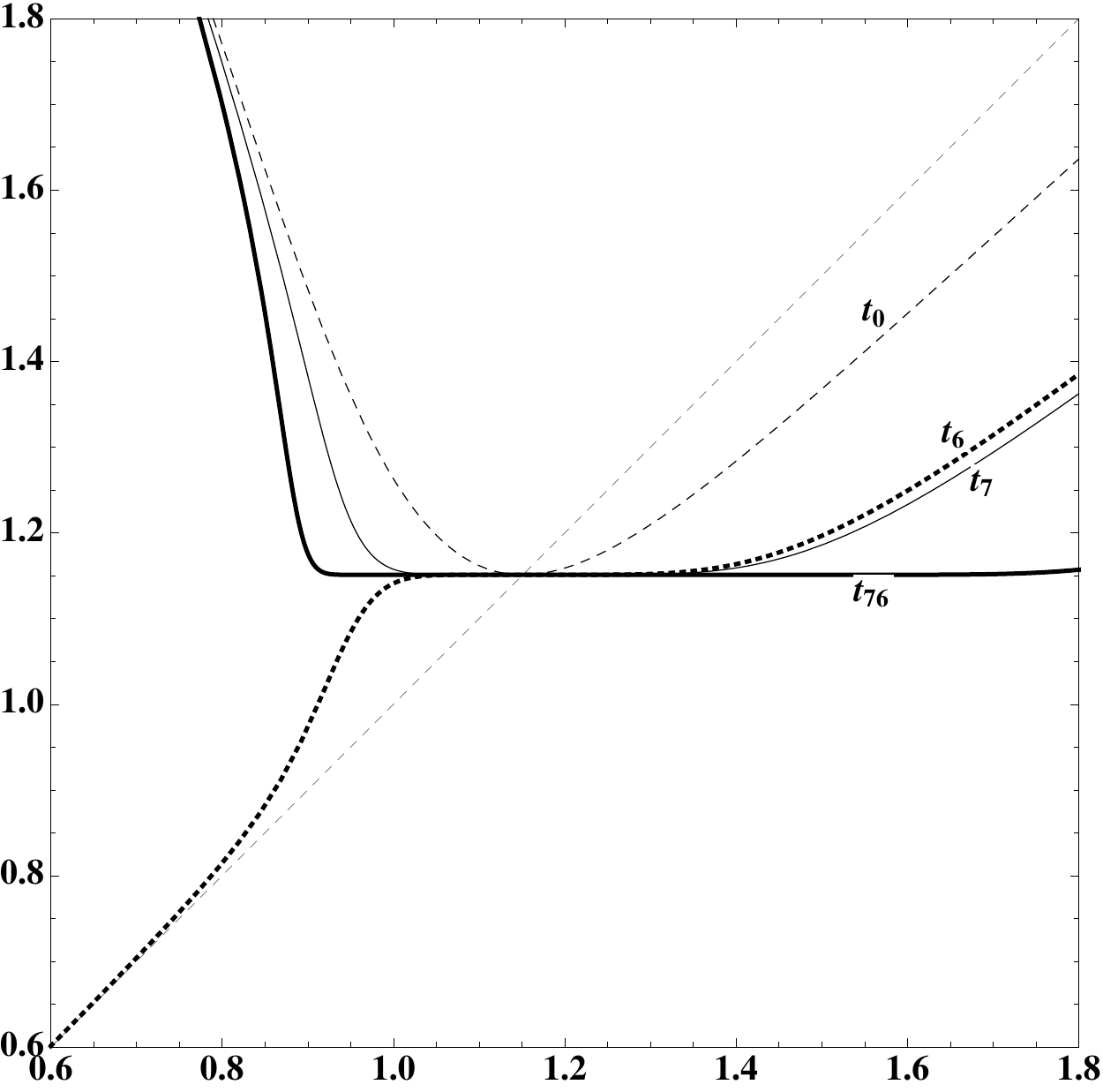}
\caption{\label{figpol2} $f(x)=x^ {11}+ 4\, x^ 2-10$. 
Comparison of results of $t_{76}$ and other iterative functions.}
\end{center}
\end{figure}

\noindent
In Fig. \ref{figpol2} we compare the graphics of the following iterative functions: $t_0$ (Newton's method),    $t_6$, $t_7$ and the composed function $t_{76}$. 
The graphic of this last function looks parallel to the $x$ axis, on a large neighborhood of $z$, which
indicates that the iterative function $t_{76}$ provides highly accurate approximations of $z$, even if we 
choose an initial approximation far from $x_0$. For example, with $x_0=2$, after 3 iterations of the Newton's method
we obtain only $0.5$ significant digits; with the same number of iterations of the $t_{76}$ iterative function
we obtain about $2410$ significant digits (see Table \ref{tabpol1}).

   \begin{table}
$$
\begin{array}{| c || c | c| c | c|  }
\hline
t_i&t_0& t_6& t_7 & t_{76} \\
\hline
s&0.5&5.3&7.6&2410.6\\
\hline
  \end{array}
  $$
  \caption{ $x_0=2$ (Example 4.4). Comparing the number of significant digits after 3 iterations.}
   \label{tabpol1}
  \end{table}
  \begin{remark}\label{obs3}
  Its is well-known that in the case of superlinear convergence the error of the $k$-th iterate $e_k=z-x_k$
  can be well approximated by the difference $x_{k+1}-x_k$. For example, the number $s$ of significant
  digits, displayed in the second row of the Table \ref{tabpol1} is in agreement with the following computations,
  when 4 iterations of $t_{76}$ are carried out, starting with  $x_0=2$: 
  $$
  \begin{array}{ll}
  x_1=t_{76} (x_0)& \Longrightarrow e_0\simeq x_1-x_0\simeq -0.799781\\
    x_2=t_{76} (x_1)& \Longrightarrow e_1\simeq x_2-x_1\simeq -0.0491500\\
      x_3=t_{76} (x_2)& \Longrightarrow e_2\simeq x_3-x_2\simeq -2.50444\times 10^ {-44}\\
         x_4=t_{76} (x_2)& \Longrightarrow e_3\simeq x_4-x_3\simeq -2.75873\times 10^ {-2411}\\
  \end{array}
  $$
  \end{remark}

\noindent
 \section{Conclusions}
 In the present article we have introduced a class of iterative methods for the numerical approximation of roots of nonlinear real functions.
 Our main goal is to propose a  recursive algorithm to construct new iterative functions $t_n$, starting with
 the classical Newton's method (to which corresponds the iterative function $t_0$), whose convergence order increases with $n$.  For each $n$, our iterative function uses a Newton-Cotes closed quadrature rule
 with $n+1$ nodes.
 We have analysed the convergence of the introduced methods, and under certain restrictions on the regularity of the considered function, we have proved that each referred method has at least convergence order $n+2$.
 The presented numerical examples illustrate the performance of the discussed methods which  can be easily extended to the case of nonlinear systems of equations. However, the analysis
 of the convergence in the  multivariate case is left for another work. We also intend in the future to explore 
 the application of the proposed methods to the solution of optimization problems.



\begin{thebibliography}{99}

\bibitem{benisrael1}{A.~Ben\--Israel},  Newton's method with modified functions, {\it  Contemporary Math. 204},  1997, 39\--50.

 \bibitem{brass} {H.~Brass and K.~Petras},   Quadrature Theory: The Theory of Numerical Integration on a Compact Interval, {\it AMS}, 2011.
 
 
\bibitem{cordero}{A.~Cordero and J.~R. Torregrosa}, Variants of Newton's method for functions of several variables, {\it  Appl. Math. Comput. 183},  2006,199\--208.

\bibitem{bjorck1} {G.~Dahlquist and \AA.~Bj\"{o}rck},   Numerical Methods in Scientific Computing, Volume I,  {\it SIAM, Philadelphia}, 2008.

\bibitem{dennis1} {J.~E. Dennis and J.~J. Mor\'{e}},   A characterization of superlinear convergence and its application to quasi\--Newton methods, {\it  Mat. Comput.}, 28, 549\--560, 1974.

 \bibitem{frontini}{M.~Frontini and E.~Sormani}, Third order methods for quadrature formulae for solving systems of nonlinear equations, {\it Appl. Math. Comput.} 149, 2004, 771\--782.


\bibitem{gautschi} {W.~Gautschi},   Numerical Analysis, An Introduction, {\it Birkh\"{a}user},  Boston, 1997.

\bibitem{isaacson} {E.~Isaacson and H.~B.~Keller},  Analysis of Numerical Methods, {\it John Wiley and Sons},  New York, 1966.

\bibitem{graca} {M.~M.~Gra\c{c}a}, Removing multiplicities in $C$ by double newtonization, {\it Appl. Math. Comput. 215(2)}, 2009, 562\--572.

\bibitem{graca1} {M.~M.~Gra\c{c}a and M.~E. ~Sousa\--Dias},   A unified framework for the computation of polynomial quadrature weights and errors, available at {\it  arXiv:1203.4795v2},  March, 2012.


\bibitem{hafiz}{M.~A. Hafiz and M.~M. Bahgat}, An efficient two\--step iterative method for solving systems of nonlinear equations,  {\it J. Math. Res.}, Vol 4., No. 4, 2012.

\bibitem{halley} {E.~Halley},   A new exact and easy method for finding the roots of equations generally and without any previous reduction, {\it Phil. Roy. Soc. London 18}, 1964, 136\--147.

\bibitem{kelley} {C.~T.~Kelley},   Iterative Methods for Linear and Nonlinear Equations, {\it SIAM}, Philadelphia, 1995.

\bibitem{krylov} {V.~I.~Krylov},   Approximate Calculation of Integrals, {\it Dover}, New York, 2005.

\bibitem{melman} {A.~Melman},   Geometry and convergence of Halley's method, {\it SIAM Rev. 39 (4)}, 1997, 728\--735.

\bibitem{mir} {N.~A.~Mir, N.~Rafiq and N.~Yasmin},   Quadrature based three\--step iterative method for nonlinear equations, {\it Gen. Math, Vol. 18, No. 4}, 2010, 31\--42.

\bibitem{rheinboldt} {W.~C. Rheinboldt},   Methods for Solving Systems of Nonlinear Equations, {SIAM, 2nd ed.}, Philadelphia,  1998.

\bibitem{thukral}{R.~Thukral}, New Sixteenth\--Order Derivative\--Free Methods for Solving Nonlinear Equations, {Amer. J. Comput. and Appl. Math. 2 (3)}, 2012, 112\--118.

\bibitem{traub} {J.~F.~Traub},   Iterative Methods for the Solution of Equations,  {\it Prentice\--Hall}, Englewood Cliffs, 1964.

\bibitem{ujevic} {N.~Ujevi\'{c}},  A method for solving nonlinear equations, {Appl. Math. Comput. 174}, 2006, 1416\--1426.


\bibitem{benisrael} {L.~ Yau and A. Ben\--Israel},   The Newton and Halley Methods for Complex Roots, {\it Amer. Math. Monthly 105}, 1998, 806\--818.
 

 \bibitem{fernando}{S. Weerakoom and T.~G. ~I. Fernando}, A Variant of Newton's Method with Accelerated Third\--Order Convergence, {\it  Appl. Math. Lett. 13},  2000, 87\--93.
 
 \bibitem{wolfram1} S. Wolfram, The Mathematica Book,  {\it Wolfram Media,  fifth ed.}, 2003.
 \end{thebibliography}
 \end{document}